\documentclass[10pt, a4paper]{article}

\usepackage{amsmath,  amsfonts}

\usepackage[T1]{fontenc}

\newtheorem{theorem}{\quad Theorem}[section]

\newcommand{\be} {\begin{equation}}
\newcommand{\ee} {\end{equation}}

\title {Uniform bound for the volume of the solutions of elliptic equations.}

\date{}

 \author{Samy Skander Bahoura\footnote {e-mail: samybahoura@gmail.com} \\ 
 {\small Equipe d'Analyse Complexe et G\'eom\'etrie.}\\  
  {\small Universit\'e Pierre et Marie Curie, 75005 Paris, France.}} 

\begin{document}

\maketitle
\begin{abstract}

We consider variational problems with regular H\"olderian weight or boundary singularity, and Dirichlet condition. We prove the boundedness of the volume of the solutions to these equations on analytic domains.

\end{abstract}



\section{Introduction and Main Results} 

We set $ \Delta = \partial_{11} + \partial_{22} $  on analytic domain $ \Omega \subset \subset {\mathbb R}^2 $.

\smallskip

We consider the following boundary value problem:

$$ (P)   \left \{ \begin {split} 
      -\Delta u & = (1+|x-x_0|^{2\beta}) V e^{u} \,\, &\text{in} \,\, & \Omega  \subset {\mathbb R}^2, \\
                  u & = 0  \,\,             & \text{on} \,\,    &\partial \Omega.              
\end {split}\right.
$$

Here:

$$ \beta \in (0, 1/2), \,\, x_0 \in \partial \Omega, $$

and,

$$ u \in W_0^{1,1}(\Omega), \,\, e^u \in L^1({\Omega}) \,\, {\rm and} \,\,  0 < a \leq V \leq b. $$

This is an equation with regular H\"olderian weight not Lipschitz in $ x_0 $ but have a weak derivative.

This problem is defined in the sense of distributions as mentioned in [10]. It arises in differents geometrical and physical situations, see for example [1, 5, 24, 28]. This type of problems was studied by many authors, with and without boundary conditions, in the subcritical case and the critical case, also for surfaces, see [1-28],  we can find in the reference existence and compactness results. In [9] we have the following important Theorem,

{\bf Theorem A}{\it (Brezis-Merle [9])}.{\it For $ u $ and $ V $ two functions relative to $ (P) $ with,
$$ 0 < a \leq V \leq b < + \infty $$
then it holds, for all compact set $ K $ of $ \Omega $:

$$ \sup_K u \leq c, $$

with $ c $ depending on $ a, b, \beta, x_0, K $ and $ \Omega $.}

We deduce from Theorem A and from the elliptic estimates that, $ u $ is uniformly bounded in $ C_{loc}^2(\Omega) $.

In [12] we have the following important Theorem and local estimate near points $ \partial \Omega \ni y\not = x_0 $ obtained by mean of the method of moving-plane for $ C^1$ functions $ V $. For $ {\cal V} $ a neighborhood of $\partial \Omega $:

{\bf Theorem B}{\it (Chen-Li [12])}.{\it For $ u $ and $ V $ two functions relative to $ (P) $ with $ V $ a $ C^1 $ function satisfying,

$$ 0 < a \leq V \leq b < + \infty, \,\, ||\nabla V||_{L^{\infty}}\leq A $$

then it holds, for all compact set $ K $ of $ {\cal V}-\{x_0\} $:

$$ \sup_K u \leq c', $$

with $ c' $ depending on $ a, b, A, \beta, x_0, K, {\cal V} $ and $ \Omega $.}

We deduce from Theorem B and from the elliptic estimates that, $ u $ is uniformly bounded in $ C_{loc}^2(\bar \Omega-\{x_0\}) $.

\smallskip

In this paper we try to prove that we have on all $ \Omega $ the boundedness of the volume of the solutions of the boundary value problem $ (P) $ if we assume that $ V $ is uniformly Lipschitz and $ C^1$ with the weight $ (1+|x-x_0|^{2\beta}), \beta \in (0,1/2) $.

Our tools are the previous two theorems A and B and a local conformal map near $ x_0 $ of the analytic domain $ \Omega $ and a Pohozaev type identity.

Here we have:

\begin{theorem} Assume that $ u $ is a solution of $ (P) $ relative to $ V $ with the following conditions:

$$  x_0 \in \partial \Omega, \,\, \beta \in (0,1/2), $$

and $ V $ a $ C^1 $ function satisfying,

$$ 0 < a \leq V \leq b,\,\, ||\nabla V||_{L^{\infty}} \leq A, $$
we have,

$$  \int_{\Omega} e^u \leq C=C(a, b, A, \beta, x_0, \Omega). $$

\end{theorem} 

\smallskip

A consequence of this theorem is a compactness result of the solutions to this Liouville type equation on analytic domains, see [4].

We have the same result if we consider the following boundary value problem on an anlaytic domain $\Omega \subset \subset {\mathbb R}^2 $:

$$ (P_{\beta})   \left \{ \begin {split} 
      -\Delta u & = |x-x_0|^{2\beta} V e^{u} \,\, &\text{in} \,\, & \Omega  \subset {\mathbb R}^2, \\
                  u & = 0  \,\,             & \text{on} \,\,    &\partial \Omega.              
\end {split}\right.
$$

Here:

$$ \beta \in (-1/2, +\infty), \,\, x_0 \in \partial \Omega, $$

and,

$$ u \in W_0^{1,1}(\Omega), \,\, |x-x_0|^{2\beta} e^u \in L^1({\Omega}) \,\, {\rm and} \,\,  0 < a \leq V \leq b. $$

Theorems A and B are true for the boundary value problem $ (P_{\beta})$ in a neighbohood of all $ y\in \bar \Omega-\{x_0\} $, because their proofs are local (see [9] and [12]) and $ x_0\in \partial \Omega $.

In this paper we try to prove that we have on all $ \Omega $ the boundedness of the volume of the solutions of the boundary value problem $ (P_{\beta}) $ if we assume that $ V $ is uniformly Lipschitz and $ C^1 $ and with the weight $ |x-x_0|^{2\beta}, \beta \in (-1/2, +\infty) $.

Our tools are the previous two theorems A and B and a local conformal map near $ x_0 $ of the analytic domain $ \Omega $ and a Pohozaev type identity.

Here we have:

\begin{theorem} Assume that $ u $ is a solution of $ (P_{\beta}) $ relative to $ V $ with the following conditions:

$$  x_0 \in \partial \Omega, \,\, \beta \in (-1/2, +\infty), $$

and $ V $ a $ C^1 $ function satisfying,

$$ 0 < a \leq V \leq b,\,\, ||\nabla V||_{L^{\infty}} \leq A, $$
we have,

$$  \int_{\Omega} |x-x_0|^{2\beta} e^u \leq C=C(a, b, A, \beta, x_0, \Omega). $$

\end{theorem} 

\smallskip

A consequence of this theorem is a compactness result of the solutions to this Liouville type equation on analytic domains, see [3].

\smallskip

{\bf Remarks:}

1) We can consider many singularities for the previous theorem 1.1, and take $ V=W\cdot \prod_{k=1}^m (1+a_k|x-x_k|^{2\beta}), x_k\in \partial \Omega,  |a_k|<<1, \beta\in (0,1/2), m \in {\mathbb N}^*, W \in C^1 $, $ W $ bounded in $ C^1 $.

2) We can consider any number of singularities for the previous theorem 1.2. $ \prod_{k=1}^m |x-x_k|^{2\beta}, x_k\in \partial \Omega, \beta \in (-1/2,+\infty) $.

3) We can consider any analytic domain, a disk or an annulus or other analytic domain.

\section{Proof of the Theorems:}

\smallskip

{\it Proof of the theorems:}

\smallskip

The proofs of theorems 1.1 and 1.2 are similar, we do the proof of theorem 1.1.

\smallskip

By corollary 1 of the paper of Brezis-Merle, we have: $ e^{ku} \in L^1(\Omega) $ for all $ k >2 $ and the elliptic estimates and the Sobolev embedding imply that: $ u \in W^{2,k}(\Omega) \cap C^{1,\epsilon}(\bar \Omega), \epsilon >0 $. By the maximum principle $ u\geq 0 $.

\smallskip

Step 1: {\it Brezis-Merle interior estimates }. By using the first eigenvalue and the first eigenfunction, with Dirichlet boundary condition, the volume is locally uniformly bounded, and thus the solutions are locally uniformly bounded by Brezis-Merle result. The solutions $ u >0 $ are locally uniformly bounded in $ C^{2,\epsilon}(\Omega) $ for $ \epsilon $ small.

\smallskip

Step 2: {\it Chen-Li boundary estimates}. Let's consider $ y\in\partial \Omega, y\not =x_0 $. Applying the moving-plane method around $ y $ and the result of Chen-Li, the solutions are uniformly bounded in a neighborhood of $ y \not =x_0 $ in $ C^{2,\epsilon}$ for $ \epsilon $ small.

\smallskip

Step 3:{\it estimates around the singularity $ x_0 $.}

\smallskip

We use a conformal map $ f=f_{x_0} $ around $ x_0 $, which maps a neighborhood of $ x_0 $, $ D_{x_0} $, to a unit half disk centered in $ 0 $, $ B^+(0, 1)$ with $ f(x_0)=0 $ and $ f(\partial \Omega \cap D_{x_0})=\{ z \in B^+(0, 1), z_1 =0 \}$. The solution, $ v=uof^{-1}$ is bounded uniformly outside $ 0 $ in $ C^2 $ norm, and is a solution of:

$$ -\Delta v= |g'(z)|^2 (1+|g(z)-x_0|)^{2\beta} W e^v, \,\, {\rm in } \,\, B^+(0,1/2), $$

with $ g=f^{-1} $ and $ W=Vog$, and $ |g'(0)|\not = 0 $.

We use a Pohozaev type identity. We use the fact that $ v $ is uniformly bounded outside $ 0 $. We multiply the equation by $ z\cdot \nabla v $ and we integrate by parts:

1) We have on a small half ball $ B^+(0,\epsilon) $:

$$ \int_{B^+(0,\epsilon)} (\Delta v)(z\cdot \nabla v) dz=\int_{B^+(0,\epsilon)} -|g'(z)|^2(1+|g(z)-x_0|^{2\beta}) W[z\cdot \nabla (e^v)]dz, $$

Thus,

$$
 \int_{\partial B^+(0,\epsilon)} (z\cdot \nabla v ) (\nabla v \cdot \nu)-\dfrac{1}{2}(z\cdot \nu ) |\nabla v|^2= $$
 
$$ = \int_{B^+(0,\epsilon)} (2+2(\beta + 1+O(\epsilon))|z|^{2\beta}) We^v dz + $$

$$ +\int_{B^+(0,\epsilon)} (1+|g(z)-x_0|^{2\beta}) [(z\cdot \nabla W)|g'(z)|^2 +(z\cdot \nabla |g'(z)|^2) W ]e^v dz+  $$

$$ - \int_{\partial B^+(0,\epsilon)} (1+|g(z)-x_0|^{2\beta}) (z\cdot \nu ) W e^v d\sigma $$

We can write, ($ v= 0 $ and $ z\cdot \nu =0 $ on $ \{z \in B^+(0,\epsilon), z_1=0\} $) and $ v $ is uniformly bounded outsode $ 0 $:

$$ \int_{\{z_1=0 \}} \frac{1}{2} (z\cdot \nu ) (\partial_{\nu} v)^2 d\sigma  + O(1)=0+O(1)= $$

$$ =\int_{B^+(0,\epsilon)} (2+2(\beta + 1+O(\epsilon))|g(z)-x_0|^{2\beta}) We^v dz + $$

$$ + \int_{B^+(0,\epsilon)} (1+|g(z)-x_0|^{2\beta}) [(z\cdot \nabla W) |g'(z)|^2 +(z\cdot \nabla |g'(z)|^2) W] e^v dz+ O(1)$$

thus, for $ \epsilon $ small enough one can compare $ z\cdot \nabla W $ and $ W $; $ |z\cdot \nabla W|\leq \epsilon A \leq a \leq W $,

$$ \int_{B^+(0,\epsilon)} |g'(z)|^2 (2+2(\beta +1+O(\epsilon))|g(z)-x_0|^{2\beta}) We^v dz =O(1), $$

Thus,

$$ \int_{B^+(0,\epsilon)} |g'(z)|^2(1+|g(z)-x_0|^{2\beta}) We^v dz =O(1), $$

and we have in a neighborhood $ D'_{x_0} $ of $ x_0$:

$$ \int_{D'_{x_0}} (1+|x-x_0|^{2\beta}) Ve^u dx =O(1), $$

uniformly.

\end{document}